\documentclass[12pt]{article}
\usepackage{amsfonts,amsbsy,graphicx}

\newcommand{\ab}{\!\!\!^{^{^\circ}}}

\newcommand{\ds}{\displaystyle}

\newcommand{\BE}{\begin{equation}}
\newcommand{\EE}{\end{equation}}
\newcommand{\fns}{\footnotesize}
\newcommand{\Lim}[1]{\lower5.5pt\hbox{${{\ds\lim}\atop ^{#1}}$}}

\begin{document}
\centerline{\Large{\bf Isoperimetric regions in  $\mathbb{H}^2$}}

\centerline{\Large{\bf between parallel horocycles}}
\ \\
\centerline{M{\fns \'ARCIO} F{\fns ABIANO DA} S{\fns ILVA}}
\ \\
{\bf Abstract.} In this work  we investigate the following 
isoperimetric problem in the hyperbolic plane: to find the 
regions of prescribed area with minimal perimeter between two 
parallel horocycles. We give an explicit and detailed description of all such regions. 
\\
\\
{\bf 1. Introduction}
\\

For a Riemannian manifold $M$, the classical isoperimetric problem
consists in classifying, up to congruency by the isometry group of
$M$, the (compact) regions $\Omega \subseteq M$ enclosing a fixed
volume that have minimal boundary volume. The existence and
regularity of solutions for a large number of cases may be
guaranteed by adapting some results from the Geometric Measure
Theory (cf. \cite{Mo2}).

For example, when $M$ is the Euclidean plane $\mathbb{R}^2$, the
classical isoperimetric problem has the disk as the unique
solution. If $M$ is a hyperbolic surface, the least-perimeter
enclosures of prescribed area are described in \cite{AM} and \cite{Si}. An
interesting version of the isoperimetric problem is to study it in
a slab. Physically, it corresponds to determine the shape of a
drop trapped between two parallel planes, which was solved by
Vogel in \cite{V}. Independently, Athanassenas studied the
isoperimetric problem between parallel planes of $\mathbb{R}^3$ in
\cite{At}. If $M$ is a slab between  two parallel horospheres in
the $3$-dimensional hyperbolic space $\mathbb{H}^{3}(-1)$, the
possible isoperimetric regions were obtained in \cite{CPS}.

In this paper we will use the upper halfplane model
$\mathbb{R}^{2}_{+}$. The parallel horocycles are
represented by the horizontal straight lines of
$\mathbb{R}^{2}_{+}$.  We will present in this paper a detailed and complete
classification of the  isoperimetric solutions.

In Section 2 we give some basic definitions in the model
$\mathbb{R}^{2}_{+}$ for the hyperbolic plane like geodesics and curves of
constant geodesic curvature. We also present a more precise
formulation for the considered isoperimetric problem and get some
preliminary characterizations by adapting the results from
\cite{CPS}. We will see
 that the possible isoperimetric regions must be delimited by such
 curves and meet the horocycles perpendicularly when this
 intersection is non-empty. We introduce the so-called geodesic
 halfdisk, horocycle halfdisk and equidistant halfdisk.

In Section 3 we read off the expressions for perimeter and area for
the regions obtained in Section 2 as the possible isoperimetric
solutions.

In Section 4 we compare the perimeter of the possible isoperimetric
regions with prescribed area. In fact, we see that this is
equivalent to investigating the regions of maximal area with
prescribed perimeter.

In Section 5 we give the isoperimetric profile for the region
between two parallel horocycles in $\mathbb{R}^{2}_{+}$ and prove
the following result:

Let $c$ be a positive real constant and $\mathcal{F}_{c} =
\{(x,y) \in \mathbb{R}^{2}_{+}: 1 \leq y \leq c \}$. Let $A>0$ and
$\mathcal{C}_{c,A}$ be the set of all  $ \Omega \subset
{\mathcal{F}_{c}}$ with area $|\Omega|=A$ and perimeter $|\partial (\Omega
\cap {{\mathcal{F\ab}}_{c}})| < \infty $, where we
suppose $\Omega$ to be connected, compact, 2-rectifiable in
$\mathcal{F}_{c}$, having as boundary (between the horocycles) a
simple rectifiable curve.

{\bf Theorem 1.1.} \it Let $L_{c,A} = \inf \{|\partial (\Omega \cap
{{\mathcal{F}\ab}_{c}})|: \Omega \in \mathcal{C}_{c,A} \}$.
Then
\begin{enumerate}
    \item there exists $\Omega \in \mathcal{C}_{c,A}$
    such that $|\partial(\Omega \cap {{\mathcal{F}\ab}_{c}})| =
    L_{c,A}$;
    \item if $\Omega \subset {\mathcal{F}_{c}}$ has minimal perimeter, the boundary of $\Omega$ has a single
    connected component made up with either
    \begin{enumerate}
            \item a halfdisk (geodesic, horocycle, equidistant) above $\{y=1\}$;
            \item a section of $\mathcal{F}_{c}$, namely
          \begin{displaymath}
           S_{[x_{0},x_{1}]}= [x_{0}, x_{1}] \times [1,c].
          \end{displaymath}
        \end{enumerate}
\end{enumerate}
More precisely, if $d$ is the hyperbolic distance between the
horocycles, we have:
\begin{description}
    \item[1.] if $d < 1$,~ there exists $A_{0}(c)$ such that
             \begin{itemize}
                \item if $A < A_{0}(c)$ then $\Omega$ is a geodesic halfdisk;
                \item if $A = A_{0}(c)$ then $\Omega$ is a geodesic halfdisk or a section;
                \item if $A > A_{0}(c)$ then $\Omega$ is a section;
              \end{itemize}
    \item[2.] if $d =1$,~ there exists  $A_{0}(c)$ such that
             \begin{itemize}
                \item if $A < A_{0}(c)$ then  $\Omega$ is a geodesic halfdisk;
                \item if $A = A_{0}(c)$ then $\Omega$ is a horocycle halfdisk or a section;
                \item if $A > A_{0}(c)$ then $\Omega$ is a section;
             \end{itemize}
    \item[3.] if $d > 1$,~ there exist two constants  $A_{0}(c) < A_{1}(c)$
    such that
             \begin{itemize}
                \item if $A < A_{0}(c)$ then $\Omega$ is a geodesic halfdisk;
                \item if $A = A_{0}(c)$ then $\Omega$ is a horocycle halfdisk;
                \item if $A_{0}(c) < A < A_{1}(c)$ then $\Omega$ is an equidistant halfdisk;
                \item if $A = A_{1}(c)$ then $\Omega$ is an equidistant halfdisk or a section;
                \item if $A > A_{1}(c)$ then $\Omega$ is a
                section.

             \end{itemize}
\end{description}
\rm

\noindent \textbf{Acknowledgments}. The present work was supported by CAPES and CNPq. The author is specially grateful to Professor Doctor Luiz Amancio Machado de Sousa Junior, from Universidade do Rio de Janeiro, for his helpful suggestions.
\\
\\
{\bf 2. Preliminaries}
\\

In this section we will introduce some basic facts and notations
that will be used along the paper. There is a large literature about the subject (we suggest  beginning with  \cite{ET}). We also adapt some important results of  \cite{CPS} to get the possible isoperimetric regions in the hyperbolic plane.

Let $\mathcal{L}^3=(\mathbb{R}^3,g)$ be the $3$-dimensional Lorentz
space endowed with the metric $g(x,y)=\ x_1 y_1 + x_2 y_2 - x_3
y_3$ and the hyperbolic plane
\begin{displaymath}
\mathbb{H}^2:=\{p=(x_1,x_2,x_3)\in \mathcal{L}^3 : \
g(p,p)=-1,\ x_3>0\}.
\end{displaymath}
We use the upper halfplane model $\mathbb{R}^{2}_{+}:= \{(x,y)
\in \mathbb{R}^{2}; y>0 \}$ for $\mathbb{H}^2$, endowed
with the metric
$<,>=ds^{2}=\displaystyle\frac{dx^{2}+dy^{2}}{y^{2}}.$

The Euclidean straight line $\{y=0\}$ is the infinity boundary of $\mathbb{R}^{2}_{+}$, denoted by $\partial_{\infty}\mathbb{R}^{2}_{+}$.

The curves of constant geodesic curvature $k \geq 0 $ in  $\mathbb{R}^{2}_{+}$ are described as follows:

\begin{enumerate}
\item {\it Geodesic:} ($k=0$). Represented by vertical Euclidean straight lines contained in  $\mathbb{R}^{2}_{+}$ and Euclidean semicircles perpendicular to
$\partial_{\infty}\mathbb{R}^{2}_{+}$  and contained in
$\mathbb{R}^{2}_{+}$;

\item {\it Geodesic circles:} ( $k>1$). Represented by Euclidean circles entirely contained in  $\mathbb{R}^{2}_{+}$;

\item {\it Horocycles:} ($k=1$). Represented by horizontal Euclidean straight lines of $\mathbb{R}^{2}_{+}$ and Euclidean circles of $\mathbb{R}^{2}_{+}$ tangent to
$\partial_{\infty}\mathbb{R}^{2}_{+}$.

\item {\it Equidistant curves:} ($0 < k < 1$). Represented by the intersection of $\mathbb{R}^{2}_{+}$ with the straight lines of $\mathbb{R}^{2}$ that are neither parallel nor perpendicular to $\{y=0\}$, and by the Euclidean circles not entirely contained in $\mathbb{R}^{2}_{+}$ and are neither tangent nor perpendicular to $\{y=0\}$.
\end{enumerate}

The isometries of $\mathbb{R}^{2}_{+}$ are the M\"obius transformations of $\widehat{\mathbb{C}}$ that leave $\mathbb{R}^{2}_{+}$ invariant. For our purposes, we are interested in the following Euclidean applications: horizontal translations, reflections with respect to a vertical geodesic, homotheties and inversions (with respect to circles centered in $\{y=0\}$).

For $\mathbb{R}^{2}_{+}$ the isoperimetric problem may be formulated as follows: ``to minimize the perimeter of a region inside two parallel horocycles (represented by two horizontal Euclidean straight lines), with prescribed area, but not counting its part of the boundary contained in the horocycles". By the {\it perimeter} of a region we mean the length of its boundary.

Since the Euclidean homothety is an isometry of $\mathbb{R}^{2}_{+}$, we take the lower horocycle as  $\{y=1\}$ to study the isoperimetric problem, so that any solution is  obtained by homothety.

By adapting the demonstration of Theorem 1.1 from \cite{CPS} to our case, namely $\mathbb{R}^{2}_{+}$, together with Lemma 2.1 of \cite {AM},  we have that there exists regular isoperimetric solutions and they are  regions whose boundary consists of curves of constant geodesic curvature  perpendicular to the horocycles (when the intersection is non-empty). Essentially, this proves the first item of our Theorem 1.1 in this present paper, stated at the Introduction.

Before we start to calculate the expressions for the perimeter and area of the regions delimited by curves of constant geodesic curvature, we present the polar coordinate system for $\mathbb{R}^{2}_{+}$ and conclude this section by giving a more precise formulation for the isoperimetric problem.

If  $(x,y)$ are the cartesian coordinates in  $\mathbb{R}^{2}_{+}$ and
 $\gamma$ is the geodesic $y>0$, we define the polar coordinates $(\rho, \theta)$ of a point $p \in \mathbb{R}^{2}_{+}$ as follows: $\rho$ is the hyperbolic distance from $p$ to  the origin $O = (0,1)$  and $\theta$ is the angle between a fixed geodesic radius $\gamma^{+}$,  given by  $\{x = 0 ;~ y \geq 1\}$, and the geodesic through $O$ and  $p$, measured counterclockwise.

The relation between these systems of coordinates is:
\begin{equation}\label{cartespolar}
   (x,y)= ~\displaystyle\frac{1}{\cosh \rho - \sinh \rho \cos \theta} ~ (\sinh \rho \sin \theta, 1),
\end{equation}
and the metric of $\mathbb{R}^{2}_{+}$ in polar coordinates is $ d\sigma^{2}=d\rho^{2}+\sinh^{2}\rho ~
    d\theta^{2}$.

We now obtain the expression for the arclength of a geodesic circle and the area of a sector as functions of the central angle $\beta$. For the sake of simplicity we take the circle of hyperbolic radius $\rho$ centered in $O$.

The geodesic circle can be parametrized by $\alpha (\theta) = (\rho, \theta)$, with constant  $\rho$ and  $0 \leq \theta \leq \beta$. Then $
d\sigma^{2}(\alpha')=\sinh^{2}\rho$. Therefore, the arclength corresponding to $\beta$ in the hyperbolic metric is
\begin{equation}\label{comprimentoesfgeod}
    L(\alpha)=\displaystyle\int_{0}^{\beta} ~
    \sqrt{d\sigma^{2}(\alpha')} ~ d\theta = \beta ~ \sinh \rho,
\end{equation}
and the area  $A$ of a sector of the disk corresponding to $\beta$ is
\begin{equation}\label{areaesfgeod}
   A=\displaystyle\int_{0}^{\beta} \!\! \displaystyle\int_{0}^{\rho} ~ \sinh \rho ~ d\rho ~ d\theta = \beta ~ (\cosh \rho -1).
\end{equation}

\begin{figure}[ht]
\begin{center}
\includegraphics[width=12cm]{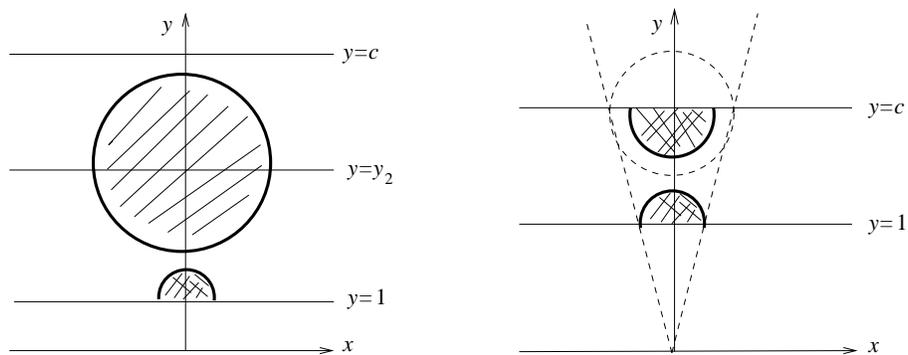}\caption{Cases 1 (left) and 2 (right).}
\label{Fig:comparadiscodisco2}
\end{center}
\end{figure}

As we mentioned above, the isoperimetric solutions are regions delimited by curves of constant geodesic curvature perpendicular to the horocycles (when the intersection is non-empty). So we have the following possibilities for barriers: vertical geodesics, geodesic circles, horocycles represented by Euclidean circles of $\mathbb{R}^{2}_{+}$ tangent to $\partial_{\infty}\mathbb{R}^{2}_{+}$, and equidistant curves represented by Euclidean circles not entirely contained in  $\mathbb{R}^{2}_{+}$ and neither tangent nor perpendicular to $\{y=0\}$. The region in ${\mathcal F}_{c}$ delimited two vertical geodesics will be called a {\it {section}}. The region in ${\mathcal F}_{c}$ delimited by geodesic circles perpendicular to $\{y=1\}$ or $\{y=c\}$ will be called {\it geodesic halfdisk}. The region in ${\mathcal F}_{c}$ delimited by horocycles and equidistant curves perpendicular to $\{y=1\}$  will be called  {\it horocycle halfdisk}  and {\it equidistant halfdisk}, respectively.  We mean by  {\it halfdisk above} (respectively {\it below}) $\{y=c\}$, the part of the Euclidean halfdisk above (respectively below) the horocycle $\{y=c\}$  (see Figure \ref{Fig:comparadiscodisco2}).

\noindent \textbf{Isoperimetric problem for
${\mathcal F}_{c}$}: \emph{fix an area value and study
the domains $\Omega\subset {\mathcal F}_{c}$ with the prescribed area
which have minimal free boundary perimeter.}

{\bf Definition 2.1:}  A (compact) minimizing region $\Omega$ for this problem will be
called an {\it isoperimetric solution} or {\it region} in
${\mathcal F}_{c}$.
\eject
{\bf 3. Expression for perimeter and area}
\\

In this section we get expressions for the perimeter and area of the possible isoperimetric solutions  $\Omega$ contained in $\mathcal{F}_c$. For our purposes we consider only the regions that are 2(-dimensional)-rectifiable (with respect to Hausdorff's measure) with boundary  1(-dimensional)-rectifiable. We denote this measure by $|\cdot|$, so that any $\Omega$ has area $|\Omega|$ and perimeter $|\partial \Omega|$, but it {\it never} counts $\partial \Omega \cap \partial \mathcal{F}_c$. For more details, see \cite{Mo2}.
\\

{\bf 3.1. Perimeter and area of a section}
\\

Let $c > 1$ and  $x_{0} < x_{1}$ be real constants. For the sake of simplicity, we consider the vertical geodesics $\{x=x_{0}\}$ and
$\{x=x_{1}\}$, contained in  $\mathbb{R}^{2}_{+}$, and the parallel horocycles $\{y=1\}$ and  $\{y=c\}$.

{\bf Lemma 3.1.1.} \it Under the notations above, if $T$ is a section then
\begin{displaymath}
      |\partial T| = 2 \ln c\,\,\,\,\,and\,\,\,\,\,|T|= (x_{1}-x_{0})(-1/c+1).
      \end{displaymath}
\rm

{\bf Proof}:~ Since the length of a vertical geodesic segment  $1 < y < c$ is
$\ln (c/1) =  \ln c$, then $|\partial T| = 2 \ln c$. And
\begin{displaymath}
|T| =  \displaystyle\int_{x_{0}}^{x_{1}} \!\!
\displaystyle\int_{1}^{c}\displaystyle\frac{1}{y^{2}} ~ dy ~ dx =
(x_{1}-x_{0})(-1/c+1).
\end{displaymath}
\ \hfill q.e.d.
\eject

{\bf 3.2. Perimeter and area for a geodesic halfdisk and a horocycle halfdisk}
\\

Consider $c \in \mathbb{R}^{*}_{+}$ and  $\{y=c\}$ a horocycle in $\mathbb{R}^{2}_{+}$. We take the Euclidean circle $S$ centered in  $(0,c)$ with radius  $r < c$. The circle  $S$ can be viewed as a geodesic circle $S_{H}$ with hyperbolic center $C_{H}=(0,h)$ and hyperbolic radius $\rho$. We want to relate the centers and the radii of $S$ and $S_{H}$.

In the hyperbolic metric, since  $C_{H}$ equidists from both  $(0,c-r)$ and $(0,c+r)$, it is easy to see that $h = \sqrt{c^{2}-r^{2}}$, $C_{H}=(0,\sqrt{c^{2}-r^{2}})$ and
\begin{equation}\label{raios}
    \rho = \displaystyle\int_{c-r}^{h}  \displaystyle\frac{1}{t} ~ dt =
    \ln \displaystyle\frac{h}{c-r}=\displaystyle\frac{1}{2} ~
    \ln \Big( \displaystyle\frac{c+r}{c-r}\Big).
\end{equation}
From (\ref{raios}) we have that $\displaystyle\frac{r}{c} = \displaystyle\frac{e^{2 \rho} -1}{e^{2
    \rho}
    +1}= \tanh \rho$.

Later we will use the relation between $|S^+|$ and $|S^-|$, where $S^+$ and $S^-$ are  halfdisks above and below  $\{y=c\}$, respectively. They are given by
\begin{displaymath}
   |S^+| = 2 \displaystyle\int_{0}^{r} \!\! \displaystyle\int_{c}^{c+\sqrt{r^{2}-x^{2}}}
   \displaystyle\frac{1}{y^{2}} ~ dy ~ dx =
   \displaystyle\frac{2}{c}
   \displaystyle\int_{0}^{r}
   \displaystyle\frac{\sqrt{r^{2}-x^{2}}}{c+\sqrt{r^{2}-x^{2}}} ~
   dx,
\end{displaymath}
and
\begin{displaymath}\label{A-}
   |S^-| = 2 \displaystyle\int_{0}^{r} \!\! \displaystyle\int_{c-\sqrt{r^{2}-x^{2}}}^{c}
   \displaystyle\frac{1}{y^{2}} ~ dy ~ dx =
   \displaystyle\frac{2}{c}
   \displaystyle\int_{0}^{r}
   \displaystyle\frac{\sqrt{r^{2}-x^{2}}}{c-\sqrt{r^{2}-x^{2}}} ~
   dx.
\end{displaymath}
Notice that $|S^-| >  |S^{+}|$. Similarly, one has 
\begin{equation}\label{comprimentos}
|\partial S^{-}| > |\partial S^{+}|.
\end{equation}

\begin{figure}[ht]
\begin{center}
\includegraphics[width=10cm]{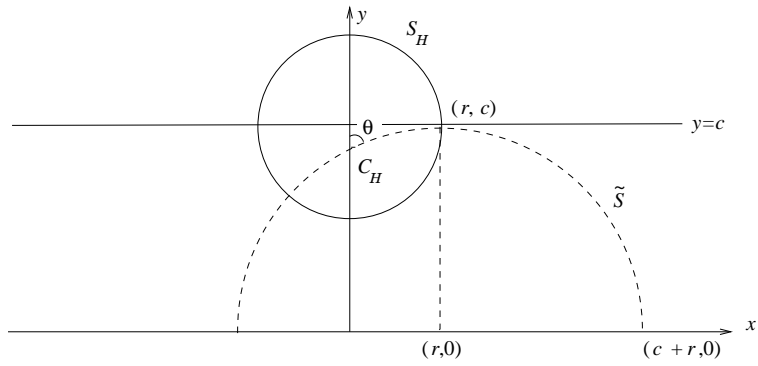}
\caption{Arc of geodesic circle corresponding to a central angle  $\theta$.} \label{Fig:esfehoro2}
\end{center}
\end{figure}

In Figure \ref{Fig:esfehoro2}, $\theta$ denotes the angle  between  $\{x=0, y \geq
h\}$ and the geodesic $\tilde{S}$ through $C_{H}$ and  $(r,c)$. Thus,  $\theta$ measures the half of the central angle corresponding to the arc of the  geodesic semicircle above $\{y=c\}$, and  $\tilde{S}$ has center  $(r,0)$ and radius $c$.


Since the Euclidean and hyperbolic metrics are conformal, in order to measure $\theta$ we  parametrize  $\tilde{S}$  as
$\alpha(t)=(c \sin t + r, c \cos t)$, with
$-\pi/2  < t < \pi/2$. Then
$C_{H}=(0,\sqrt{c^{2}-r^{2}}) = \alpha(t_{0})=(c
\sin t_{0} + r, c \cos t_{0})$ so that $\sin t_{0}= -r/c$ and
\begin{equation}\label{costheta}
   \cos \theta = - \sin t_{0} = r/c.
\end{equation}

If $\bar{S}$ is the region delimited by $\tilde{S}$, axis $y$ and $\{y=c\}$ then
\begin{displaymath}
    |\bar{S}| = \displaystyle\int_{h}^{c}
    \displaystyle\frac{r-\sqrt{c^{2}-y^{2}}}{y^{2}} ~ dy =
    -r/c + \pi/2 - \arcsin(h/c).
\end{displaymath}

\begin{figure}[ht]
\begin{center}
\includegraphics[width=11cm]{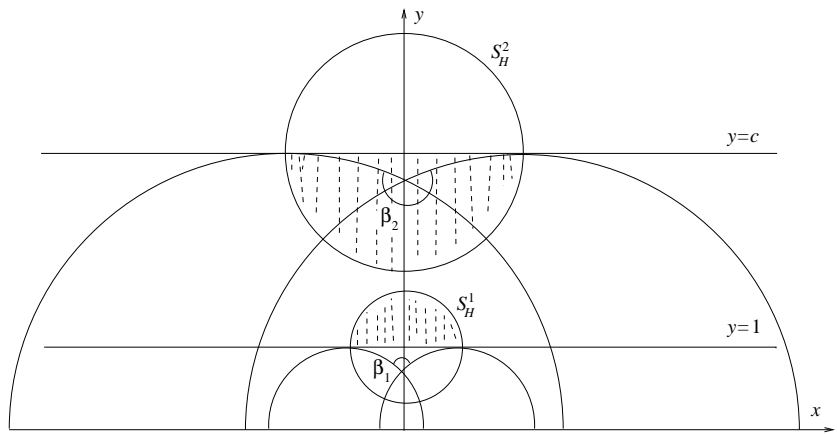}
\caption{Perimeter and area for geodesic halfdisks.}
\label{Fig:comparacompesfe}
\end{center}
\end{figure}

Suppose  $c > 1$ and consider the parallel horocycles $\{y=1\}$ and $\{y=c\}$. Let $S_{1}$ be the circle centered at  $(0,1)$ with radius $r_{1} < 1$, and $S_{2}$ the circle  centered at  $(0,c)$ with radius  $r_{2} <
c-1$ (see Figure \ref{Fig:comparacompesfe}). Hence,  $S_{1}$ can be viewed as a geodesic circle $S_{H}^{1}$ with hyperbolic center $(0,h_{1})=(0,\sqrt{1-r_{1}^{2}})$ and hyperbolic radius $ \rho_{1} = \displaystyle\frac{1}{2} ~
    \ln \Big( \displaystyle\frac{1+r_{1}}{1-r_{1}}\Big)$, and  $S_{2}$ as a geodesic circle  $S_{H}^{2}$ with hyperbolic center $(0,h_{2})=(0,\sqrt{c^{2}-r_{2}^{2}})$ and hyperbolic radius $ \rho_{2} = \displaystyle\frac{1}{2} ~
    \ln \Big( \displaystyle\frac{c+r_{2}}{c-r_{2}}\Big)$.


Let  $\beta_{1}$ be the central angle of $S_{H}^{1}$ corresponding  to the  arc above  $\{y=1\}$, and  $\beta_{2}$ the central angle of $S_{H}^{2}$ corresponding  to the arc  below $\{y=c\}$. By (\ref{costheta}) we have $\beta_{1} = 2 \arccos (r_{1})$ and $\beta_{2} = 2 \pi - 2 \arccos(r_{2}/c).$

For geodesic halfdisks it holds the following result (see Figure \ref{Fig:comparacompesfe}):

{\bf Lemma 3.2.1.} \it Under the notations above, let $\tilde{S_{1}}$ be the geodesic through $C_{H}^{1}=(0,h_{1})$ and  $(r_{1},1)$, and
$\tilde{S_{2}}$ the geodesic through
$C_{H}^{2}=(0,h_{2})$ and  $(r_{2},c)$. Let
$\theta_{1}= \beta_{1}/2$ and  $\theta_{2}=
\pi - \beta_{2}/2$,~  $0 < \theta_{1},
\theta_{2} < \pi/2$. Let $S_{1}^{+}$ be the geodesic halfdisk delimited by $S_{H}^{1}$ and above $\{y=1\}$, and $S_{2}^{-}$ the halfdisk delimited by  $S_{H}^{2}$ and below  $\{y=c\}$. Then
\begin{equation}\label{L1L2}
      |\partial S_{1}^{+}|= 2 \theta_{1} \cot \theta_{1},\,\,\,\,\,
      |\partial S_{2}^{-}|= 2 (\pi - \theta_{2}) ~ \cot \theta_{2},
   \end{equation}
and
\begin{equation}\label{areashiprb}
       |S_{1}^{+}| =  \displaystyle\frac{ 2 \theta_{1}}{\sin \theta_{1}} -\pi + 2 \cos \theta_{1}, \,\,\,\,\,
        |S_{2}^{-}| =   \displaystyle\frac{2 (\pi -  \theta_{2})}{\sin \theta_{2}} -\pi - 2 \cos \theta_{2}.
    \end{equation}
\rm

{\bf{Proof}}:~ By (\ref{comprimentoesfgeod}), the arclengths determined by $\beta_{1}$ and  $\beta_{2}$ are
\begin{displaymath}
    |\partial S_{1}^{+}|=  \beta_{1} \sinh \rho_{1}\,\,\,\,\,{\rm and}\,\,\,\,\,
    |\partial S_{2}^{-}|=  \beta_{2} \sinh
    \rho_{2}.
    \end{displaymath}

We have 
\begin{displaymath}
    \begin{array}{ll}
         \sinh \rho_{1} = \sinh \Big( \displaystyle\frac{1}{2} ~
    \ln \Big( \displaystyle\frac{1+r_{1}}{1-r_{1}}\Big)\Big) =
    \displaystyle\frac{r_{1}}{\sqrt{1-r_{1}^{2}}}, \\
    \sinh \rho_{2} = \sinh \Big(\displaystyle\frac{1}{2} ~
    \ln \Big( \displaystyle\frac{c+r_{2}}{c-r_{2}}\Big)\Big) =
    \displaystyle\frac{r_{2}}{\sqrt{c^{2}-r_{2}^{2}}}.
    \end{array}
\end{displaymath}

Since $\cot \theta_{1} = \displaystyle\frac{r_{1}}{\sqrt{1-r_{1}^{2}}}$ and $\cot \theta_{2} =
      \displaystyle\frac{r_{2}}{\sqrt{c^{2}-r_{2}^{2}}}$, the first part of the lemma is proved.

Now  observe that $|S_{1}^{+}|/2 = |\tilde{S_{1}}| - |\bar{S}_{1}|$, where $\tilde{S_{1}}$ is the sector corresponding to
$\theta_{1}$ and  $\bar{S_{1}}$ is the region delimited by
$\tilde{S_{1}}$, axis $y$ and the horocycle $\{y=1\}$. In the same way, $|S_{2}^{-}|/2 = |\tilde{S_{2}}| + |\bar{S}_{2}|$,
where $\tilde{S_{2}}$ is the sector corresponding to
$ \beta_{2}/2 = \pi  - \theta_{2}$ and
$\bar{S}_{2}$ is the region delimited by
$\tilde{S_{2}}$, axis  $y$ and the horocycle $\{y=c\}$.

Therefore, by (\ref{areaesfgeod})
\begin{equation}\label{comparaarea}
    \begin{array}{ll}
    |S_{1}^{+}| =  2 \theta_{1} ~(\cosh \rho_{1} -1) -2 (-r_{1} + \pi/2 -\arcsin
    (h_{1})),\\
    |S_{2}^{-}|=  2 (\pi -  \theta_{2}) ~(\cosh \rho_{2}
      -1) + 2( -r_{2}/c + \pi/2-\arcsin(h_{2}/c)).
    \end{array}
\end{equation}

But
\begin{equation}\label{cosh}
    \begin{array}{ll}
         \cosh \rho_{1} = \cosh \Big(\displaystyle\frac{1}{2} ~
    \ln \Big( \displaystyle\frac{1+r_{1}}{1-r_{1}}\Big)\Big) =
    \displaystyle\frac{1}{\sqrt{1-r_{1}^{2}}}, \\
    \cosh \rho_{2} = \cosh \Big(\displaystyle\frac{1}{2} ~
    \ln \Big( \displaystyle\frac{c+r_{2}}{c-r_{2}}\Big)\Big) =
    \displaystyle\frac{c}{\sqrt{c^{2}-r_{2}^{2}}}.
    \end{array}
\end{equation}

Since $r_{1}^{2}+ \Big(\sqrt{1-r_{1}^{2}}\Big)^{2} = 1$ and $(r_{2}/c)^{2}+
(\sqrt{c^{2}-r_{2}^{2}}/c)^{2}=1$,
we have
\begin{equation}\label{arcos}
    \begin{array}{ll}
        \arccos(r_{1})= \arcsin ( \sqrt{1-r_{1}^{2}}) = \arcsin
    (h_{1}),\\
        \arccos(r_{2}/c)= \arcsin \Big(
        \displaystyle\frac{\sqrt{c^{2}-r_{2}^{2}}}{c}\Big) = \arcsin
    (h_{2}/c).
    \end{array}
\end{equation}

Furthermore, by (\ref{costheta}) it follows that $\cos \theta_{1} = r_{1}$ and $\cos \theta_{2} = r_{2}/c$.

Therefore,
\begin{equation}\label{seno}
       \sin \theta_{1} = \sqrt{1-r_{1}^{2}}\,\,\,\,\,{\rm and}\,\,\,\,\,
        \sin \theta_{2} = \displaystyle\frac{\sqrt{c^{2}-r_{2}^{2}}}{c}.
   \end{equation}

By (\ref{comparaarea}), (\ref{cosh}), (\ref{arcos}) and
(\ref{seno}), the proof of  (\ref{areashiprb}) is complete\hfill q.e.d.
\\

We observe that a horocycle $H$ can be viewed as a limit geodesic circle with hyperbolic center in $\partial_{\infty}\mathbb{R}^{2}_{+}$. By
(\ref{costheta}), we have $\cos \theta_{1} = r_{1}$ and the horocycle is obtained when  $r_{1}$ converges to 1, that is, $\theta_{1}$ converges to 0. Hence we get the expressions for the perimeter and the area of the horocycle halfdisk as the following consequence of Lemma 3.2.1:

{\bf Corollary 3.2.2.} \it  Let  $H$ be the horocycle halfdisk above $\{y=1\}$ represented by a Euclidean semicircle with center  $(0,1)$ and radius $1$. Then
\begin{equation}\label{compareahoro}
   |\partial H|= 2\,\,\,\,\,and\,\,\,\,\,|H| = 4 - \pi.
   \end{equation}
\rm

{\bf{Proof}}:~ It is enough to calculate $|\partial S_{1}^{+}|$ and $|S_{1}^{+}|$
from (\ref{L1L2}) and (\ref{areashiprb}) for the limit case when
$\theta_{1} \rightarrow 0$\hfill q.e.d.
\ \\

{\bf 3.3. Perimeter and area for an equidistant halfdisk}
\\

 Let $\bar{E}$ be the equidistant curve represented by a Euclidean circle with center  $(0,1)$ and radius  $r > 1$. The Euclidean equation of $\bar{E}$ is given by $x^{2}+(y-1)^{2}=r^{2}$. Then $\bar{E} \cap
\partial_{\infty}\mathbb{R}^{2}_{+} =
\{(-\sqrt{r^{2}-1},0),(\sqrt{r^{2}-1},0)\}$. The curve $\bar{E}$ is equidistant from the geodesic $\eta$ with equation $x^{2}+y^{2}=r^{2}-1$.
If $\rho$ denotes the hyperbolic distance between  $\bar{E}$ and
$\eta$, then $\rho$ is the hyperbolic  distance between $(0,1+r)$ and $(0,\sqrt{r^{2}-1})$, so that $\rho = \ln \Big( \displaystyle\frac{r+1}{r-1}\Big)^{\frac{1}{2}}$, whence 
\begin{equation}\label{rhoequi}
    r = \coth \rho.
\end{equation}
If  $\alpha$ is the non-oriented angle between  $\bar{E}$ and
$\eta$,~ $0 < \alpha < \pi/2$, then (for instance, see
Proposition 3 in Chapter 5 of \cite{ET})
\begin{equation}\label{rhoalfa}
    \tanh \rho = \sin \alpha.
\end{equation}

{\bf Lemma 3.3.1.} \it  Under the notations above, let $E$ be the equidistant halfdisk above $\{y=1\}$. Then
\begin{equation}\label{compareaequid}
   \begin{array}{ll}
   |\partial E|= \displaystyle\frac{2}{\cos \alpha} ~ \ln\Big(\displaystyle\frac{1}{\sin \alpha} + \cot \alpha\Big),\\
   ~\\
   |E| = \displaystyle\frac{2}{\sin \alpha}  - \pi + \displaystyle\frac{2}{\cot \alpha}~
   \ln\Big(\displaystyle\frac{1}{\sin \alpha} + \cot \alpha\Big).
   \end{array}
\end{equation}
\rm

{\bf{Proof}:}~ In order to calculate $|\partial E|$, we parametrize
$E$ by
\begin{displaymath}
\beta(t) = (r \cos t, 1+ r \sin t), ~ 0 \leq t \leq \pi.
\end{displaymath}

Then
\begin{displaymath}
\begin{array}{ll}
|\partial E|= 2
\displaystyle\int_{0}^{\pi/2}
\displaystyle\frac{r}{1+ r \sin t} ~ dt= \displaystyle\frac{2 r}{\sqrt{r^{2}-1}} ~ \ln
\Big(\displaystyle\frac{r+\tan (t/2)-\sqrt{r^{2}-1}}{r+\tan (t/2)+\sqrt{r^{2}-1}}\Big)\Bigg|_{0}^{\pi/2}\\
~\\
=\displaystyle\frac{2 r}{ \sqrt{r^{2}-1}} ~ \ln (r
+ \sqrt{r^{2}-1}).
\end{array}
\end{displaymath}

By (\ref{rhoequi}) and (\ref{rhoalfa}) we have $r = 1/\sin \alpha$, whence $
\sqrt{r^{2}-1 }= \cot \alpha$, because $0 < \alpha < \pi/2$. Therefore, $
|\partial E| = \displaystyle\frac{2}{\cos \alpha}
 ~ \ln \Big(\displaystyle\frac{1}{\sin \alpha} + \cot \alpha
 \Big)$, and the first part of (\ref{compareaequid}) is proved. Now, 
\begin{displaymath}
|E| = 2 \displaystyle\int_{0}^{r} \!\!
\displaystyle\int_{1}^{1+\sqrt{r^{2}-x^{2}}}
   \displaystyle\frac{1}{y^{2}} ~ dy ~ dx =  2r - \pi + \displaystyle\frac{1}{\sqrt{r^{2}-1}}
    \ln \Big| \displaystyle\frac{r \sqrt{r^{2}-1} + (r^{2}-1)}{r \sqrt{r^{2}-1} -
    (r^{2}-1)}\Big|.
\end{displaymath}

By (\ref{rhoequi}) and (\ref{rhoalfa}), it follows that $|E|$, as function of the equidistance angle $\alpha$, is given by $|E| = \displaystyle\frac{2}{\sin \alpha}  - \pi +
\displaystyle\frac{2}{\cot \alpha}~
   \ln\Big(\displaystyle\frac{1}{\sin \alpha} + \cot \alpha\Big)$, which proves the second part of (\ref{compareaequid})\hfill q.e.d.

\ \\
{\bf 4. Comparison of perimeters of regions with prescribed area}
\\

In this section we analyze the expressions of the perimeter and area of the regions delimited by curves of constant geodesic curvature. Their isoperimetric profiles in  $\mathcal{F}_{c}$ will be obtained in the next section as  functions of its hyperbolic width $d$. Since we have been taking the horocycles $\{y=1\}$ and $\{y=c\}$, the constant $c$ must satisfy the following condition: if $H$ is a horocycle halfdisk above $\{y=c\}$ and $T$ is a section in $\mathcal{F}_{c}$, then $|\partial H|=|\partial T|$. By (\ref{compareahoro}), this means $2 = 2 \ln c$, whence $c = e$ and $d=1$. This is why we compare $d$ with $1$ in Theorem 1.1.

Let $S_{1}$ be a Euclidean circle with radius $r_{1}$ and center
$(0,1)$ above $\{y=1\}$. When $0 < r_{1} < 1$, ~ $S_{1}$
delimits a geodesic halfdisk. Consider the limit cases $\theta_{1}
\rightarrow 0$ and $\theta_{1} \rightarrow \pi/2$, which correspond to a
horocycle halfdisk ($r_{1} \rightarrow 1$) and a point ($r_{1}
\rightarrow 0$), respectively. From  (\ref{costheta}) and (\ref{L1L2})  we have
\begin{equation}\label{limcomp1}
\Lim{\theta_{1} \to  0}{|\partial S_{1}^{+}|} = 2,\,\,\,\,\,
\Lim{\theta_{1} \to \pi/2}{|\partial S_{1}^{+}|} = 0,
\end{equation}
and from (\ref{cartespolar}),
\begin{equation}\label{limarea1}
\Lim{\theta_{1} \to 0}{|S_{1}^{+}|} = 4 - \pi,\,\,\,\,\,
\Lim{\theta_{1} \to \pi/2}{|S_{1}^{+}|} = 0.
\end{equation}

Therefore, $|\partial S_{1}^{+}|$ increases from 0 to 2  when $r_{1}$ varies from 0 to 1, while $|S_{1}^{+}|$ increases from  0 to $4 - \pi$.

If $r_{1} > 1$ we have an equidistant halfdisk. By (\ref{rhoequi}) and
 (\ref{rhoalfa}) the limit cases are obtained for  $\alpha \rightarrow
\pi/2$ (that is,  $r_{1}\rightarrow 1$)  and $\alpha \rightarrow 0$ (that is, $r_{1}\rightarrow \infty$). By (\ref{compareaequid}), 
\begin{equation}\label{limcompequid2}
\Lim{\alpha \to \pi/2}{|\partial E|} = 2,\,\,\,\,\,
\Lim{\alpha \to 0}{|\partial E|} =
\infty,
\end{equation}
and
\begin{equation}\label{limareaequid2}
\Lim{\alpha \to \pi/2}{|E|} = 4 -\pi,\,\,\,\,\,
\Lim{\alpha \to 0}{|E|} =
\infty.
\end{equation}

Therefore, both $|\partial E|$ and $|E|$ increase infinitely while $r_{1}$ increases.

Since  $c > 1$, let  $S_{2}$ be a Euclidean semicircle with radius $r_{2}$ and center $(0,c)$ below  $\{y=c\}$. When  $0 < r_{2} < c$, ~ $S_{2}$ delimits a geodesic halfdisk. By (\ref{costheta}), the limit cases  
$\theta_{2} \rightarrow 0$ and $\theta_{2} \rightarrow \pi/2$   correspond to $r_{2} \rightarrow c$ and 
 $r_{2}
\rightarrow 0$, respectively. By (\ref{L1L2}), we have for these limit cases
\begin{equation}\label{limcomp2}
\Lim{\theta_{2} \to 0}{|\partial S_{2}^{-}|} = \infty,\,\,\,\,\,
\Lim{\theta_{2} \to \pi/2}{|\partial S_{2}^{-}|} = 0,
\end{equation}
and by (\ref{areashiprb}), 
\begin{equation}\label{limarea2}
\Lim{\theta_{2} \to 0}{|S_{2}^{-}|} = \infty,\,\,\,\,\,
\Lim{\theta_{2} \to \pi/2}{|S_{2}^{-}|} =
0.
\end{equation}

Then we see that $|\partial S_{2}^{-}|$ and $|S_{2}^{-}|$ increase infinitely while $r_{2} \rightarrow c$. If  $r_{2} \geq c$,~ $S_{2}$ delimits a horocycle halfdisk or an equidistant halfdisk. By
(\ref{limcomp2}) and (\ref{limarea2}), if  $r_{2} \geq c$,
then both $|\partial S_{2}^{-}|$ and $|S_{2}^{-}|$ diverge to infinity.

From the analysis we have just done for the perimeter and area of the possible isoperimetric solutions, there are only the following cases to consider:
\begin{enumerate}
    \item to compare a geodesic halfdisk above  $\{y=1\}$ with a geodesic disk entirely contained in $\mathcal{F}_c$;
    \item to compare a geodesic halfdisk above  $\{y=1\}$ with a  geodesic halfdisk below $\{y=c\}$;
    \item  to compare a horocycle halfdisk above  $\{y=1\}$ with a  geodesic halfdisk below $\{y=c\}$;
    \item  to compare an equidistant halfdisk above  $\{y=1\}$ with a  geodesic halfdisk below $\{y=c\}$.
\end{enumerate}

In order to prove the second part of Theorem 1.1, one must determine the least-perimeter regions with prescribed area. For this purpose, we will use a strategy: we  determine the regions with prescribed perimeter and biggest area. In fact, it is enough to show that if a region has the maximum  area among all regions with a prescribed perimeter, then it has the minimum perimeter among all regions with the same prescribed area (see Lemma 4.1 below). Since we have just listed  all  possible isoperimetric solutions besides the section, Lemma 4.1 will then  refer to the above case 2. The other cases are proved analogously. Without loss of generality we suppose that the geodesic  halfdisk above $\{y=1\}$ has maximum area when compared to any geodesic halfdisk below  $\{y=c\}$  with the same perimeter.

{\bf Lemma 4.1.} \it Let  $\Omega_{0}$ be the geodesic halfdisk above $\{y=1\}$ with $|\Omega_{0}| \geq |\Omega|$, whenever $|\partial \Omega|= |\partial \Omega_{0}|$, for any geodesic halfdisk $\Omega$ below  $\{y=c\}$,~ $c > 1$. If $\Omega_{1}$ is a geodesic halfdisk below 
$\{y=c\}$ with  $|\Omega_{0}| = |\Omega_{1}|$, then 
$|\partial \Omega_{0}| \leq |\partial\Omega_{1}|$.\rm

{\bf Proof:}~ Suppose by contradiction  that $|\partial \Omega_{0}| >
|\partial \Omega_{1}|$. By (\ref{limcomp2}) we can increase the radius of the Euclidean circle that represents  $\Omega_{1}$ till we get a geodesic halfdisk $\Omega'$ such that $|\partial \Omega'|= |\partial \Omega_{0}|$. This procedure could fail if $\Omega'$ surpassed $\{y=1\}$, but then the section will prevail as the isoperimetric solution. This fact will be proved later on in Section 5. By (\ref{limarea2}), the area increases with the radius. Therefore, $|\Omega'| >
|\Omega_{1}| = |\Omega_{0}|$ and $|\partial \Omega'|= |\partial \Omega_{0}|$. This  is a contradiction with the fact that  $\Omega_{0}$ maximizes the area when compared to regions of the same perimeter, by hypothesis \hfill q.e.d.

Till the end of this section we are going to compare the area of the possible isoperimetric solutions for a  prescribed perimeter.

For case 1 described above, we compare the area of a geodesic halfdisk above  $\{y=1\}$  with a geodesic disk entirely contained in $\mathcal{F}_c$, when they have the same perimeter. Let $S$ be the Euclidean circle with radius  $r_{2}$, ~ $0 < r_{2} <
    y_{2}-1$, and center $(0, y_{2})$,~ $1 < y_{2} <  c$, which delimits the geodesic halfdisk (see Figure \ref{Fig:comparadiscodisco2}).


Consider  $\theta_{2} ~ \in ~ ]0, \pi/2[$ such that
$\cos \theta_{2} = r_{2}/c$. By
(\ref{comprimentoesfgeod}), (\ref{areaesfgeod}) and (\ref{raios}), if $\mathcal{S}$ is the geodesic disk corresponding to a central angle of  $2 \pi$  then $|\partial \mathcal{S}|= 2 \pi  \cot \theta_{2}$ and  $|\mathcal{S}| =
\displaystyle\frac{2 \pi}{\sin \theta_{2}} -2 \pi$.

By (\ref{L1L2}), (\ref{areashiprb}) and the information from the previous paragraph, we show that $ |S_{1}^{+}|>|\mathcal{S}|$  when  $|\partial S_{1}^{+}|=|\partial \mathcal{S}|$ in the next Lemma.

{\bf Lemma 4.2.} \it Let $\theta_{1},\theta_{2} ~ \in ~ ]0, \pi/2[$ such that
\begin{equation}\label{xcotan0}
    \theta_{1} \cot \theta_{1} = \pi  \cot \theta_{2}.
\end{equation}
Then
\begin{equation}\label{xsin0}
  \displaystyle\frac{2 \theta_{1}}{\sin \theta_{1}} + 2 \cos \theta_{1}  - \pi >
   \displaystyle\frac{2 \pi}{\sin \theta_{2}} -2 \pi.
\end{equation}
\rm

{\bf {Proof}:}~ For $\theta_{2} ~ \in ~ ]0,
\pi/2[$, by calculating the squares of
(\ref{xcotan0}) and using that  $\cos^{2} \theta_{2} = 1
- \sin^{2} \theta_{2}$, one has $1/\sin \theta_{2} = \sqrt{\theta_{1}^{2} \cot^{2} \theta_{1} +
    \pi^{2}}/\pi$. Thus
\begin{equation}\label{termodearea}
    \displaystyle\frac{2 \pi}{\sin \theta_{2}} -2 \pi = 2 \sqrt{\Big(\displaystyle\frac{\theta_{1}}{\sin
    \theta_{1}}\Big)^{2} - \theta_{1}^{2} + \pi^{2}} -
    2 \pi.
\end{equation}

Now we replace the right-hand side  of  (\ref{xsin0}) by (\ref{termodearea}), and define 
\begin{displaymath}
 A (\theta_{1} ):= \displaystyle\frac{2 \theta_{1}}{\sin \theta_{1}} + 2 \cos \theta_{1}  - \pi -
 2 \sqrt{\Big(\displaystyle\frac{\theta_{1}}{\sin
    \theta_{1}}\Big)^{2} - \Big(\theta_{1}\Big)^{2} + \pi^{2}} +
    2 \pi,
\end{displaymath}
so that (\ref{xsin0}) will hold if and only if $A(\theta_{1})>0$.

We observe that
\begin{equation}\label{limitesdearea}
   \Lim{\theta_{1} \to 0} {A (\theta_{1})}  = 4 + \pi - 2 \sqrt{\pi^{2} +
    1} > 0\,\,\,\,\,{\rm and}\,\,\,\,\,
    \Lim{\theta_{1} \to \pi/2}{A (\theta_{1}
    )}= 0.
   \end{equation}
Moreover,
\begin{displaymath}
    \displaystyle\frac{d A (\theta_{1} )}{d \theta_{1} }
    = \displaystyle\frac{2 \cos \theta_{1} \Big( \sin \theta_{1} \cos \theta_{1} - \theta_{1} \Big)}{\sin^{2}
    \theta_{1}}\Big\{ 1 - \displaystyle\frac{\theta_{1}}{\sqrt{\theta_{1}^{2} + \Big( \pi^{2} - \theta_{1}^{2} \Big)
    \sin^{2} \theta_{1}}} \Big\} < 0,
\end{displaymath}
because $\theta_{1} ~ \in ~ ]0, \pi/2[$ implies $\cos \theta_{1} >0 $, ~$\sin \theta_{1} \cos \theta_{1} -
\theta_{1} < 0$ and
\begin{displaymath}
 0 < \displaystyle\frac{\theta_{1}}{\sqrt{\theta_{1}^{2} + \Big( \pi^{2} - \theta_{1}^{2} \Big)
    \sin^{2} \theta_{1}}} < 1.
\end{displaymath}
Therefore, $A (\theta_{1} )$ decreases in  $]0,
\pi/2[$. By (\ref{limitesdearea}),
we conclude that $A (\theta_{1} ) > 0 $ in  $]0,
\pi/2[$, whence (\ref{xsin0}) is proved \hfill q.e.d.

We conclude from Lemma 4.2 that the geodesic halfdisk above $\{y=1\}$ is the isoperimetric solution, instead of the geodesic disk, which concludes case 1.

Now we study case 2. By (\ref{L1L2}) and (\ref{areashiprb}), we will show in the next Lemma and Corollary that $|S_{1}^{+}|>|S_{2}^{-}|$ when  $|\partial S_{1}^{+}|=|\partial S_{2}^{-}|$. In  Figure
    \ref{Fig:comparadiscodisco2}, the dashed circle was obtained from the lower by a Euclidean homothety so that the corresponding geodesic halfdisks have the same perimeter. By (\ref{comprimentos}), in order to have $|\partial S_{1}^{+}|=|\partial S_{2}^{-}|$, it is necessary to decrease the radius of $S_{2}^{-}$.


{\bf Lemma 4.3.} \it  Let $\theta_{1},\theta_{2} ~ \in ~ ]0,\pi/2]$ such that 
\begin{equation}\label{xcotan}
    \theta_{1} \cot \theta_{1} = (\pi - \theta_{2}) \cot \theta_{2}.
\end{equation}
Then
\begin{equation}\label{xsin}
  \displaystyle\frac{\theta_{1}}{\sin \theta_{1}} + \cos \theta_{1} \geq  \displaystyle\frac{\pi - \theta_{2}}
  {\sin \theta_{2}} - \cos \theta_{2}.
\end{equation}
\rm

{\bf {Proof}:}~For $\theta_{1},\theta_{2} \in  ]0,
\pi/2]$, we define $f(\theta_{1}) =
\displaystyle\frac{\theta_{1}}{\sin \theta_{1}} + \cos \theta_{1}$
and  $F(\theta_{1},\theta_{2}) = f(\theta_{1}) + f(\theta_{2}) -
\displaystyle\frac{\pi}{\sin \theta_{2}}$. We want to show that $F(\theta_{1},\theta_{2}) \geq 0$. By  (\ref{xcotan}) we can define $\theta_{1}$ implicitly as a function of $\theta_{2}$. Namely, we get a function  $g$ such that $\theta_{1} = g(\theta_{2})$. Let $h_{1}(\theta_{2}) = F(g(\theta_{2}),\theta_{2})$ and
$h_{2}(\theta_{2}) = h_{1}(\theta_{2}) \sin \theta_{2} + \pi =
\sin \theta_{2} ~ f(g(\theta_{2})) + \sin \theta_{2} ~
f(\theta_{2})$.

The function  $h_{2}(\theta_{2})$ is  $\mathcal{C}^{\infty}$ and
\begin{equation}\label{um}
    h_{2}'(\theta_{2}) = \cos \theta_{2} ~ f(\theta_{1}) + \sin \theta_{2} ~ f'(\theta_{1}) ~ g'(\theta_{2}) + \cos \theta_{2} ~ f(\theta_{2}) + \sin \theta_{2} ~
    f'(\theta_{2}).
\end{equation}
Hence
\begin{equation}\label{dois}
   f'(\theta_{1}) = \displaystyle\frac{\sin \theta_{1} - \theta_{1} \cos \theta_{1}}{\sin^{2} \theta_{1}} -
   \sin \theta_{1} = \displaystyle\frac{- \cos \theta_{1} (2\theta_{1} - \sin 2\theta_{1})}{2 \sin^{2}
   \theta_{1}}.
\end{equation}
From the Implicit Function Theorem we have
\begin{displaymath}
   g'(\theta_{2}) = \displaystyle\frac{- \cot \theta_{2} - (\pi -\theta_{2}) \csc^{2}\theta_{2}}{\cot \theta_{1}  - \theta_{1} \csc^{2}
   \theta_{1}}.
\end{displaymath}
Now observe that $\cot \theta_{1}  - \theta_{1} \csc^{2}
   \theta_{1} = \displaystyle\frac{\sin 2\theta_{1}  - 2\theta_{1}}{2
   \sin^{2}\theta_{1}}$.
Therefore,
\begin{equation}\label{quatro}
   g'(\theta_{2}) = \displaystyle\frac{2 (\pi -\theta_{2}) + \sin 2\theta_{2}}{2 \sin^{2}\theta_{2}} ~ \displaystyle\frac{2 \sin^{2}\theta_{1}}{2\theta_{1} - \sin
   2\theta_{1}}.
\end{equation}
By substituting (\ref{dois}) and (\ref{quatro}) in  (\ref{um}), we obtain
\begin{equation}\label{cinco}
   h_{2}'(\theta_{2}) = 2 \cos^{2} \theta_{2} + \displaystyle\frac{\theta_{1} \cos \theta_{2}}{\sin \theta_{1}}-  \displaystyle\frac{(\pi - \theta_{2}) \cos \theta_{1}}{\sin
   \theta_{2}}.
\end{equation}
Since  $ h_{2}'(\theta_{2}) =  h_{1}'(\theta_{2}) \sin \theta_{2} +
h_{1}(\theta_{2}) \cos \theta_{2}$, it follows from  (\ref{cinco}) that
\begin{equation}\label{seis}
\begin{array}{lll}
  h_{1}'(\theta_{2}) \sin \theta_{2} = 2 \cos^{2} \theta_{2} + \displaystyle\frac{\theta_{1} \cos \theta_{2}}{\sin
  \theta_{1}}- \displaystyle\frac{(\pi -\theta_{2}) \cos \theta_{1}}{\sin \theta_{2}} - F(\theta_{1},\theta_{2}) \cos \theta_{2}
  \\
  ~\\
  = \cos^{2} \theta_{2} - \cos \theta_{2} ~ \cos \theta_{1} - \displaystyle\frac{(\pi -\theta_{2})}{\sin
  \theta_{2}}(\cos \theta_{1} - \cos \theta_{2})\\
  ~\\
  = (\cos \theta_{2} - \cos \theta_{1}) (\cos \theta_{2} + \displaystyle\frac{\pi -\theta_{2}}{\sin
  \theta_{2}}).
   \end{array}
\end{equation}
For $\theta_{1},\theta_{2} ~ \in ~ ]0,
\pi/2]$, we have $ \theta_{1} \leq \pi -
\theta_{2}$, which by (\ref{xcotan}) implies
\begin{equation}\label{sete}
   \cot \theta_{1}= \Big(\displaystyle\frac{\pi -\theta_{2}}{\theta_{1}}\Big) \cot \theta_{2} \geq \cot \theta_{2} \Rightarrow \theta_{1} \leq \theta_{2} \Rightarrow \cos \theta_{1} \geq \cos \theta_{2}.
\end{equation}
By (\ref{seis})  and (\ref{sete}) we have $h_{1}'(\theta_{2})
\leq 0$. Thus $F(\theta_{1},\theta_{2}) = F(g(\theta_{2}),\theta_{2}) =$ $h_{1}(\theta_{2})$ is a decreasing function for $\theta_{2} ~ \in ~ ]0,
\pi/2]$. By  (\ref{xcotan}), for $\theta_{2}
= \pi/2$ and $\theta_{1} = g
(\pi/2)$, one has $\cos (g
(\pi/2)) = 0$ and therefore $g
(\pi/2) = \pi/2$. Since
\begin{displaymath}
h_{1}(\pi/2) =
F(g(\pi/2),\pi/2) =
F(\pi/2,\pi/2)=0,
\end{displaymath}
then $F(\theta_{1},\theta_{2})\geq
F(g(\pi/2),\pi/2) = 0$. Consequently, $f(\theta_{1}) + f(\theta_{2}) \geq \displaystyle\frac{\pi}{\sin
\theta_{2}}$, whence (\ref{xsin}) is proved.

The equality occurs if and only if
$\theta_{1}=\theta_{2}=\pi/2$. 

\hfill q.e.d.

{\bf Corollary 4.4.} \it Let  $\theta_{1},\theta_{2} ~ \in ~ ]0,
\pi/2[$ such that
\begin{equation}\label{xcotan2}
    \theta_{1} \cot \theta_{1} = (\pi - \theta_{2}) \cot \theta_{2}.
\end{equation}
Then
\begin{equation}\label{xsin2}
  \displaystyle\frac{\theta_{1}}{\sin \theta_{1}} + \cos \theta_{1} >  \displaystyle\frac{\pi - \theta_{2}}
  {\sin \theta_{2}} - \cos \theta_{2}.
\end{equation}
\rm

We finally conclude from Corollary 4.4 that the geodesic halfdisk above $\{y=1\}$ is the isoperimetric solution for case 2.

\begin{figure}[ht]
\begin{center}
\includegraphics[width=14cm]{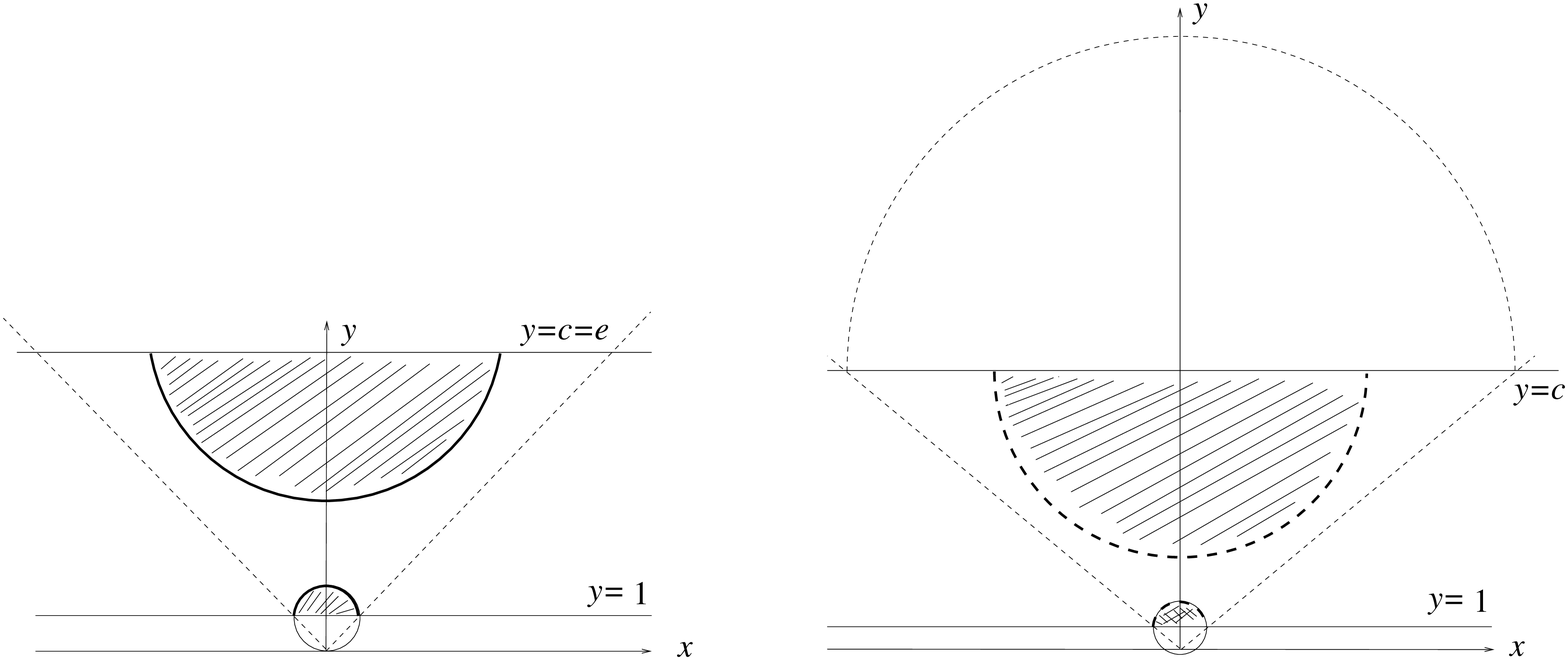} \caption{Cases 3 (left) and 4 (right).}
\label{Fig:comparahorodisco}
\end{center}
\end{figure}

By (\ref{compareahoro}), (\ref{L1L2}) and 
    (\ref{areashiprb}), we show in the next Lemma that $|H| > |S_{2}^{-}|$  when $|\partial H| = |\partial S_{2}^{-}|$. Case 3 is illustrated in Figure  \ref{Fig:comparahorodisco}.

{\bf Lemma 4.5.} \it  Let $\theta_{2} ~ \in ~ ]0, \pi/2[$
such that $1 = (\pi - \theta_{2}) \cot \theta_{2}$. Then $2 >  \displaystyle\frac{\pi - \theta_{2}}{\sin \theta_{2}} - \cos \theta_{2}$.
\rm

{\bf {Proof}:}~ Since the horocycle is obtained  from the geodesic halfdisk above $\{y=1\}$ when  $\theta_{1}
\rightarrow 0$, then it is enough to make $\theta_{1} \rightarrow 0$ in (\ref{xcotan2}) and (\ref{xsin2}). The result follows from the continuity of the involved functions \hfill q.e.d.

We conclude from Lemma 4.5 that the horocycle halfdisk above $\{y=1\}$ is the isoperimetric solution, instead of the geodesic  halfdisk below $\{y=c\}$.

Now we analyze case 4. By (\ref{compareaequid}), (\ref{L1L2}) and (\ref{areashiprb}), we show in the next Lemma that $|E| >|S_{2}^{-}|$  when $|\partial E| = |\partial S_{2}^{-}|$. In Figure  \ref{Fig:comparahorodisco}, the dashed circle was obtained from the lower by a Euclidean homothety so that they have the same perimeter. In order to have $|\partial E| = |\partial S_{2}^{-}|$, it is necessary to decrease the radius of $S_{2}^{-}$.


{\bf Lemma 4.6.} \it Let $\alpha, \theta_{2} ~ \in ~ ]0, \pi/2[$
such that
\begin{equation}\label{xcotan4}
    \displaystyle\frac{1}{\cos \alpha} ~\ln \Big(\displaystyle\frac{1}{\sin \alpha} + \cot \alpha \Big) = (\pi - \theta_{2}) \cot \theta_{2}.
\end{equation}
Then 
\begin{equation}\label{xsin4}
 \displaystyle\frac{1}{\sin \alpha}  + \displaystyle\frac{1}{\cot \alpha}~
   \ln\Big(\displaystyle\frac{1}{\sin \alpha} + \cot \alpha\Big) \geq  \displaystyle\frac{\pi - \theta_{2}}
   {\sin \theta_{2}} - \cos \theta_{2}.
\end{equation}
\rm

{\bf {Proof}:}~ For $\alpha,\theta_{2} ~ \in ~ ]0,
\pi/2[$  we define
\begin{displaymath}
F(\alpha,\theta_{2}) = \displaystyle\frac{1}{\sin \alpha}  +
\displaystyle\frac{1}{\cot \alpha}~
   \ln\Big(\displaystyle\frac{1}{\sin \alpha} + \cot \alpha\Big) +
   \displaystyle\frac{\theta_{2}}{\sin \theta_{2}} + \cos
   \theta_{2} - \displaystyle\frac{\pi}{\sin \theta_{2}}.
\end{displaymath}
By (\ref{xcotan4}), we can implicitly define $\theta_{2}$ as a function of  $\alpha$. Namely, one gets a function $g$ such that
$\theta_{2} = g(\alpha)$. Consider the functions $h_{1}(\alpha) = F(\alpha, g(\alpha))$, $h_{2}(\alpha)= h_{1}(\alpha) \sin \alpha$. Then \begin{displaymath} h_{2}(\alpha)= \Big(
\displaystyle\frac{1}{\sin \alpha} + \displaystyle\frac{1}{\cot
\alpha}~
   \ln\Big(\displaystyle\frac{1}{\sin \alpha} + \cot \alpha\Big)
\Big) \sin\alpha - \Big( \displaystyle\frac{\pi - g(\alpha)}{\sin
g(\alpha)} - \cos g(\alpha) \Big) ~ \sin \alpha.
\end{displaymath}

The function $h_{2}(\alpha)$ is  $\mathcal{C}^{\infty}$ and
\begin{equation}\label{doisum}
\begin{array}{cl}
    h_{2}'(\alpha) = \displaystyle\frac{\sin\alpha}{\cos^{2}\alpha}~
   \ln\Big(\displaystyle\frac{1}{\sin \alpha} + \cot \alpha\Big) - \displaystyle\frac{\sin\alpha}{\cos\alpha} +
   \sin\alpha ~ \ln\Big(\displaystyle\frac{1}{\sin \alpha} + \cot \alpha\Big) +\\
   ~\\
   + g'(\alpha)~
    \displaystyle\frac{\cos\theta_{2} (\sin(2 \theta_{2}) + 2 (\pi - \theta_{2}))}{2 \sin^{2} \theta_{2}} \sin\alpha -
    \displaystyle\frac{(\pi - \theta_{2})}{\sin\theta_{2}} \cos\alpha + \cos\theta_{2} \cos\alpha.
\end{array}
\end{equation}
 From the Implicit Function Theorem we have 
\begin{displaymath}
   g'(\alpha) = - \displaystyle\frac{\displaystyle\frac{\sin\alpha}{\cos^{2}\alpha}~
   \ln\Big(\displaystyle\frac{1}{\sin \alpha} + \cot \alpha\Big)- \displaystyle\frac{1}{\sin\alpha \cos\alpha}}
   {\cot\theta_{2} + (\pi - \theta_{2})\csc^{2} \theta_{2}}.
\end{displaymath}
Since  $\cot\theta_{2} + (\pi - \theta_{2})\csc^{2}
\theta_{2} = \displaystyle\frac{\sin 2\theta_{2}  + 2 (\pi -
\theta_{2})}{2 \sin^{2}\theta_{2}}$,  then
\begin{equation}\label{doisquatro}
   g'(\alpha) = \Big(\displaystyle\frac{2 \sin^{2}\theta_{2}}{\sin 2\theta_{2}  + 2 (\pi -
\theta_{2})}\Big) ~ \Big(\displaystyle\frac{1}{\sin\alpha
\cos\alpha}-\displaystyle\frac{\sin\alpha}{\cos^{2}\alpha}~
   \ln\Big(\displaystyle\frac{1}{\sin \alpha} + \cot \alpha\Big)\Big).
\end{equation}
By substituting (\ref{doisquatro}) in  (\ref{doisum}), we obtain
\begin{equation}\label{doiscinco}
\begin{array}{ll}
 h_{2}'(\alpha) = \displaystyle\frac{\sin\alpha}{\cos^{2}\alpha}~
   \ln\Big(\displaystyle\frac{1}{\sin \alpha} + \cot \alpha\Big) - \displaystyle\frac{\sin\alpha}{\cos\alpha} +
   \sin\alpha ~ \ln\Big(\displaystyle\frac{1}{\sin \alpha} + \cot \alpha\Big) +\\
   ~\\
   + \displaystyle\frac{\cos\theta_{2}}{\cos\alpha} -\displaystyle\frac{\sin^{2}\alpha}{\cos^{2}\alpha}
   ~\cos\theta_{2} ~
  \ln\Big(\displaystyle\frac{1}{\sin \alpha} + \cot \alpha\Big) -
    \displaystyle\frac{(\pi - \theta_{2})}{\sin\theta_{2}} \cos\alpha + \cos\theta_{2} \cos\alpha.
\end{array}
\end{equation}
Since  $ h_{2}'(\alpha) =  h_{1}'(\alpha) \sin \alpha +
h_{1}(\alpha) \cos \alpha$, it results from  (\ref{doiscinco}) that
\begin{displaymath}
  h_{1}'(\alpha) \sin\alpha = \displaystyle\frac{1 - \sin\alpha ~
  \cos\theta_{2}}{\cos\alpha} ~ \Big( \tan\alpha ~
  \ln\Big(\displaystyle\frac{1}{\sin \alpha} + \cot \alpha\Big) - \displaystyle\frac{1}{\sin\alpha}\Big).
\end{displaymath}
For  $\alpha \in ~ ]0, \pi/2[$, if $l(\alpha) = \tan\alpha ~ \ln(1/\sin
\alpha + \cot \alpha)$, then $l'(\alpha) = \sec^{2}\alpha ~ k(\alpha)$, where
$k(\alpha) = \ln(1/\sin \alpha + \cot
\alpha) - \cos\alpha$. Since $k'(\alpha)= -\cos^{2} \alpha/\sin
\alpha < 0$, then  $k(\alpha)$ is decreasing in  $]0,
\pi/2[$ and $\Lim{\alpha \to \pi/2}{k(\alpha)} = 0$. So
$k(\alpha) > 0$ in  $]0, \pi/2[$.
Consequently, $l(\alpha)$ is increasing in  $]0,
\pi/2[$. Since $\Lim{\alpha \to \pi/2}{l(\alpha)} = 1$, then
$l(\alpha) < 1$ and therefore
\begin{displaymath}
\tan\alpha ~
  \ln\Big(\displaystyle\frac{1}{\sin \alpha} + \cot \alpha\Big) - \displaystyle\frac{1}{\sin\alpha} < 0.
\end{displaymath}
Moreover, for $\alpha, \theta_{2} \in ~ ]0,
\pi/2[$  we have  $0< \sin\alpha
  \cos\theta_{2} < 1$. Thus, for $\alpha \in ~ ]0, \pi/2[$ 
we conclude that $h_{1}'(\alpha) \sin\alpha < 0$ and therefore
$h_{1}'(\alpha) < 0$. Namely, $h_{1}$ is decreasing in  $]0,
\pi/2[$. In particular,
\begin{displaymath}
h_{1}(\alpha) \geq \Lim{\alpha \to \pi/2}{ h_{1}(\alpha)} = \Lim{\alpha
\to \pi/2}{F(\alpha,g(\alpha))} = 2 -
\Big(\displaystyle\frac{\pi-\beta}{\sin\beta} - \cos\beta \Big),
\end{displaymath}
where  $\beta = \Lim{\alpha \to \pi/2}{g(\alpha)}$. From  (\ref{compareahoro})  it follows that 
\begin{displaymath}
\displaystyle\frac{|H| + \pi}{2} = 2,
\end{displaymath}
where $H$ is a horocycle halfdisk above $\{y=1\}$. By (\ref{areashiprb}),
\begin{displaymath}
\displaystyle\frac{|G| + \pi}{2} =
\displaystyle\frac{\pi-\beta}{\sin\beta} - \cos\beta,
\end{displaymath}
where  $G$ is a geodesic halfdisk with the same perimeter as $H$ (just take $\alpha \rightarrow
\pi/2$ in  (\ref{xcotan4}) and $\theta_{2}=\beta$). But from Lemma 4.5 we conclude that
\begin{displaymath}
h_{1}(\alpha) = 2 - \Big(\displaystyle\frac{\pi-\beta}{\sin\beta}
- \cos\beta \Big) > 0,
\end{displaymath}
whence (\ref{xsin4}) is proved \hfill q.e.d.

From Lemma 4.2, Corollary 4.4, Lemma 4.5 and Lemma 4.6, we conclude that the family of geodesic, horocycle and equidistant halfdisks above $\{y=1\}$ are the solutions to the isoperimetric problem, instead of the geodesic halfdisks below  $\{y=c\},~ ~c > 1$.

\ \\
{\bf 5.  Isoperimetric Profile in  $\mathbb{R}^{2}_{+}$}
\\

\begin{figure}[ht]
\begin{center}
\includegraphics[width=12cm]{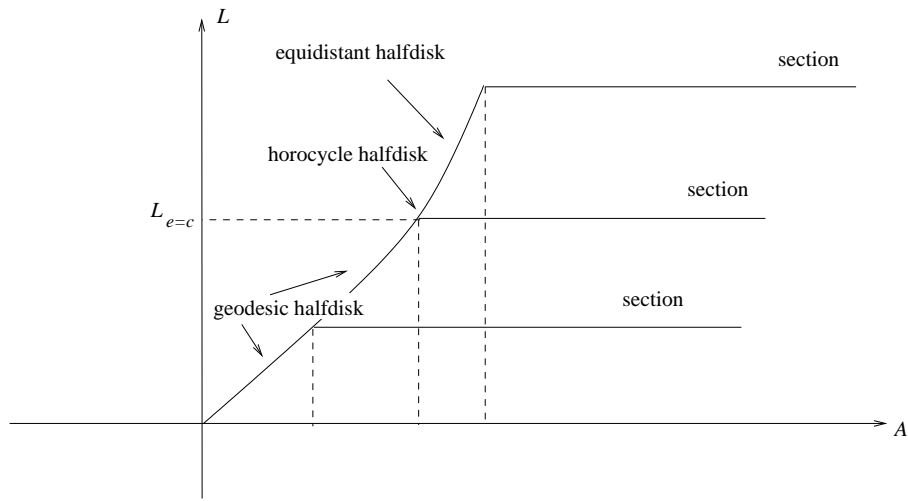}
\caption{Isoperimetric profile for the region between the parallel horocycles.} \label{Fig:perfiliso}
\end{center}
\end{figure}

In this section we study the isoperimetric profile for $\mathcal{F}_c$ (see Figure \ref{Fig:perfiliso}). We adapt a well-known result from the Isoperimetric Problem Theory which  guarantees that the boundaries of the connected components of an isoperimetric solution are curves with the same constant geodesic curvature (for instance, see Lemma 2.1 of \cite{AM}). Before showing that a minimizing region is made up with a single connected component, we prove that a connected component of an isoperimetric region {\it  must} be either a section or a halfdisk above the horocycle $\{y=1\}$. Here we need (\ref{limcomp1})-(\ref{limarea2}). The perimeter of the section in $\mathcal{F}_c$  is equal to $2 \ln c$. Now there are only three possibilities that we classify according to the hyperbolic distance $d=\ln c$.

{\bf First Possibility: $d <1$} 

\begin{figure}[ht]
\begin{center}
\includegraphics[width=12cm]{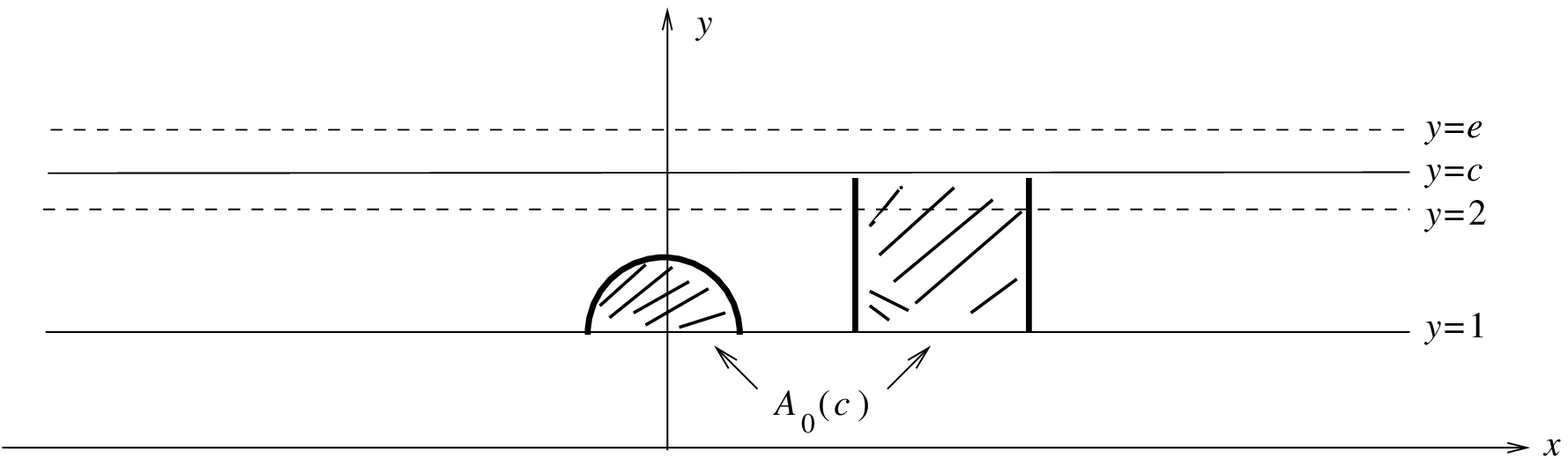} \caption{Case $c < e$.}
\label{Fig:caso1teo}
\end{center}
\end{figure}

\begin{figure}[ht]
\begin{center}
\includegraphics[width=12cm]{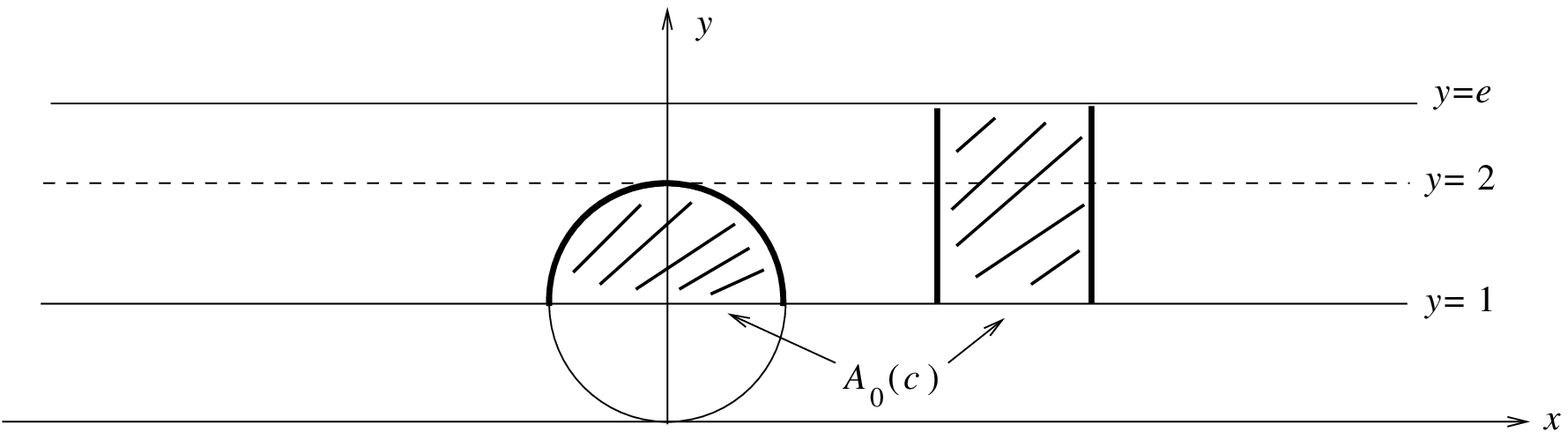} \caption{Case $c = e$.}
\label{Fig:caso2teo}
\end{center}
\end{figure}

\begin{figure}[ht]
\begin{center}
\includegraphics[width=12cm]{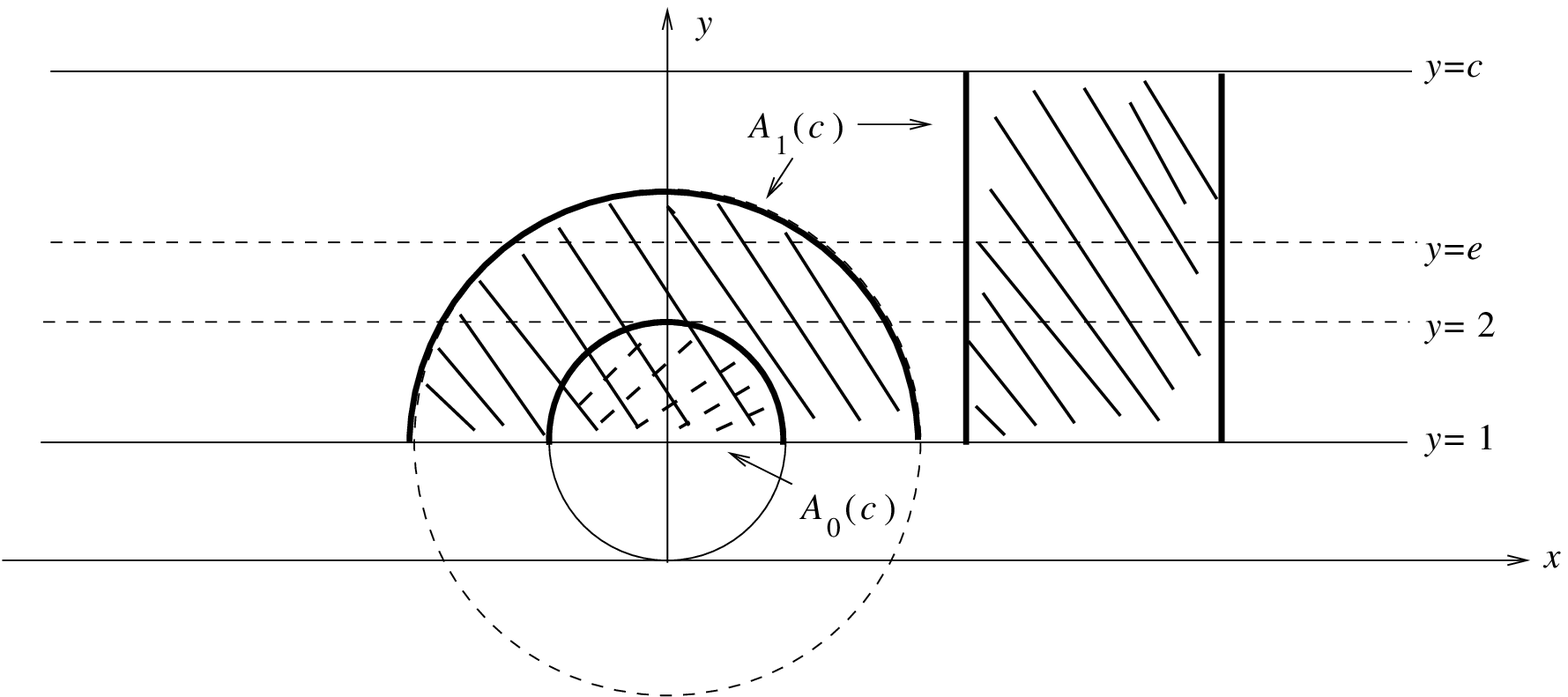} 
\caption{Case $c > e$.} \label{Fig:caso3teo}
\end{center}
\end{figure}

\begin{enumerate}
    \item Consider a horocycle $\{y=c\}$  with  $1 < c < e$. Let  $A_{0}(c)$ be the area of the geodesic halfdisk  $S_{0}$ above $\{y=1\}$, centered at  $(0,1)$ with Euclidean radius $r_{0}(c)$ and $|\partial S_0|=|\partial T_0|$, where  $T_{0}$ is a section with $|T_0|=A_{0}(c)$ (see Figure
    \ref{Fig:caso1teo}). Since  $c < e$, then $|\partial T_{0}| < 2$ (which is the perimeter of the horocycle halfdisk above  $\{y=1\}$). 


Consequently,
\begin{itemize}
    \item if $A=A_{0}(c)$ then $|\partial S_{0}|=|\partial T_{0}|$ and $|S_{0}|=|T_{0}|=A$. Therefore, the minimizing region  $\Omega$ is a geodesic halfdisk or a section;
    \item if  $A<A_{0}(c)$, let  $S_{1}$ be a geodesic halfdisk with area $A$, centered at $(0,1)$ and with Euclidean radius $r_{1}$. Since both $|S_1|$ and $|\partial S_1|$ decrease with $r_1$, we have  $r_{1} < r_{0}(c)$ and $|\partial S_{1}| < |\partial S_{0}|$.  Let $T_{1}$ be a section with $|T_1|=A$. Then  $|S_{1}|=|T_{1}|=A$, but $|\partial S_{1}| < |\partial T_{1}|= |\partial T_{0}|=|\partial S_{0}|$. Therefore, the minimizing $\Omega$ is a geodesic halfdisk. In this case, we observe that $|\Omega|=A<|S_{0}|$, so that  $\Omega$ can neither be  a horocycle nor an equidistant halfdisk;
    \item if  $A>A_{0}(c)$, let $S_{2}$ be a geodesic halfdisk with $|S_{2}|=A$, centered at  $(0,1)$ and with Euclidean radius $r_{2}$. Since both $|S_{2}|$ and $|\partial S_{2}|$ increase with $r_2$, then $r_{2} > r_{0}(c)$ and $|\partial S_{2}| > |\partial S_{0}|$. Let  $T_{2}$ be a section with $|T_{2}|=A$. Then  $|S_{2}|=|T_{2}|=A$, but $|\partial S_{2}| > |\partial T_{2}|=|\partial T_{0}|=|\partial S_{0}|$. Therefore, the minimizing $\Omega$  is a section.
    \end{itemize}

{\bf Second Possibility: $d =1$} 

    \item Suppose $d=1$. Consider the horocycle
    $\{y=c\}$ with  $c =e$. Then $A_{0}(c)=4 - \pi$ is the area of the horocycle halfdisk $S_{0}$ above $\{y=1\}$, centered at  $(0,1)$ with Euclidean radius $r_{0}(c)=1$ and $|\partial S_{0}|=|\partial T_{0}|$, where $T_0$ is a section with $|T_0|=A_{0}(c)$ (see Figure
    \ref{Fig:caso2teo}).  In this case, $|\partial T_{0}|=2$.


Consequently,
\begin{itemize}
    \item if $A=A_{0}(c)$, then  $|S_{0}|=|T_{0}|=A$. Therefore,  the minimizing $\Omega$ is a horocycle halfdisk or a section;
    \item if  $A<A_{0}(c)$, let  $S_{1}$ be a geodesic halfdisk with $|S_1|=A$,  centered at  $(0,1)$ and with Euclidean radius $r_{1}$. Since both $|S_1|$ and $|\partial S_1|$ increase with $r_1$ till it becomes a horocycle disk, then  $r_{1} < 1$  and $|\partial S_1| < |\partial S_0|$. Let  $T_{1}$ be a section with $|T_1|=A$. Then  $|S_{1}|=|T_{1}|=A$, but $|\partial S_1|<|\partial T_1|=|\partial T_0|=|\partial S_0|$. Therefore, the minimizing $\Omega$ is a geodesic halfdisk;
    \item if  $A>A_{0}(c)$, let  $S_{2}$ be an equidistant halfdisk $|S_2|=A$, centered at  $(0,1)$ and with Euclidean radius $r_{2}$. Since both $|S_2|$ and $|\partial S_2|$ increase infinitely with $r_2$, then  $r_{2} > 1$ and $|\partial S_{2}|>|\partial S_{0}|$. Let  $T_{2}$ be a section with $|T_2|=A$. Then  $|S_{2}|=|T_{2}|=A$, but $|\partial S_{2}|>|\partial T_{2}|=|\partial T_{0}|=|\partial S_{0}|$. Therefore, the minimizing  $\Omega$ is a section.
\end{itemize}

{\bf Third Possibility: $d > 1$}

    \item Suppose $d > 1$. Consider a horocycle $\{y=c\}$ with  $ c > e$. Let  $A_{0}(c)=4 - \pi$ be the area of the horocycle halfdisk  $S_{0}$ above  $\{y=1\}$, centered at  $(0,1)$ with Euclidean radius  $r_{0}(c)=1$ and $|\partial S_0|=2$. Let  $T_{0}$ be a section with $|T_0|=A_{0}(c)$ and  $A_{1}(c)$ be the area of an equidistant halfdisk  $S_{1}$ above $\{y=1\}$, centered at $(0,1)$ with Euclidean radius $r_{1}(c)$ and $|\partial S_1|=|\partial T_{1}|$, where $T_1$ is a  section with $|T_1|=A_{1}(c)$ (see Figure \ref{Fig:caso3teo}). In this case, we observe that $|\partial T_{1}|> 2$.


Consequently,
\begin{itemize}
    \item if  $A=A_{0}(c)=4 - \pi$ then $|S_{0}|=|T_{0}|=A$, but $|\partial T_{0}| > 2 =|\partial S_{0}|$. Therefore, the minimizing  $\Omega$ is a horocycle halfdisk;
    \item if  $A=A_{1}(c)$ then  $|S_{1}|=|T_{1}|=A$ and $|\partial S_{1}|=|\partial T_{1}|$. Therefore, the minimizing $\Omega$ is an equidistant halfdisk or a section;
    \item if  $A<A_{0}(c)$, let  $S_{2}$ be a geodesic halfdisk with $|S_2|=A$, centered at  $(0,1)$  and with Euclidean radius $r_{2}$. Then $r_{2}< r_{0}(c)$ and $|\partial S_{2}|<|\partial S_{0}|$. Let  $T_{2}$ be a section with $|T_2|=A$. Then $|S_{2}|=|T_{2}|=A$, but  $|\partial S_{2}| < |\partial T_{2}|=|\partial T_{0}|=|\partial S_{0}|$. Therefore, the minimizing $\Omega$ is a geodesic halfdisk;
    \item if  $A_{0}(c)<A<A_{1}(c)$, let  $S_{3}$ be an equidistant halfdisk with $|S_3|=A$, centered at  $(0,1)$ and with Euclidean radius $r_{3}$. Then  $r_{0}(c) < r_{3} < r_{1}(c)$ and  $|\partial S_{3}| < |\partial S_{1}|$. Let $T_{3}$ be a section with $|T_3|=A$. Then $|S_{3}|=|T_{3}|=A$, but $|\partial S_{3}| < |\partial T_{3}|= |\partial T_{1}|=|\partial S_{1}|$. Therefore, the minimizing  $\Omega$ is an equidistant halfdisk;
    \item if  $A>A_{1}(c)$, let  $S_{4}$ be an equidistant halfdisk with $|S_4|=A$, centered at  $(0,1)$ and with Euclidean radius  $r_{4}$. Then  $r_{4} > r_{1}(c)$ and $|\partial S_{4}|>|\partial S_{1}|$.  Let  $T_{4}$ be a section with $|T_4|=A$. Then $|S_{4}|=|T_{4}|=A$, but $|\partial S_{4}| > |\partial T_{4}|= |\partial T_{1}|=|\partial S_{1}|$. Therefore, the minimizing  $\Omega$ is a section.
\end{itemize}
\end{enumerate}

REMARK 5.1: A minimizing region consists of only one connected component, and in fact it is enough to show that it can not have two. If this were the case, their geodesic curvatures would agree. Consider  $A >0$ and  $\Omega'$ a region with area $A$ and two disjoint sections. Their ``gluing" would result in another section with area $A$ but with smaller perimeter, because two vertical geodesics would not count anymore. Then $\Omega'$ is not minimizing.

The other case to consider is two connected components consisting of two geodesic halfdisks above $\{y=1\}$. In this case, we use the fact that  a non-regular region is not minimizing: let $A
> 0$ and $\Omega'$ be a region with area  $A$ and two geodesic halfdisks above $\{y=1\}$ with the same Euclidean radius, hence the same geodesic curvature. By sliding one of them over $\{y=1\}$ till it touches the other, since  horizontal translations are isometries of the hyperbolic plane, we get a non-regular region  $\Omega''$
with area  $A$. Then  $\Omega''$ does not have the least-perimeter among all regions  with prescribed area $A$. Since $|\Omega'|=|\Omega''|$, $\Omega'$ is not  minimizing.

Therefore, a minimizing region must consist of a single connected component.

Now we prove Theorem 1.1. 

{\bf {Proof}:}~ The first part of Theorem 1.1 was already discussed in the Preliminaries. The existence of such an isoperimetric region follows from adaptions of some results from  \cite{Mo1} and  \cite{Mo2}: the group $G$ of isometries of $\mathbb{R}^{2}_{+}$ that leave $\mathcal{F}_c$ invariant consists of  horizontal Euclidean translations and Euclidean reflections with respect to a vertical geodesic, so that $\mathcal{F}_c/G$ is homeomorphic to the interval  $[0,1]$, hence compact.

The second part of Theorem 1.1 follows from the analysis of the isoperimetric profile done in the three possibilities above, together with REMARK 5.1.

\ \\
M\'arcio Fabiano da Silva\\
Universidade Federal do ABC\\
r. Catequese 242, 3rd floor\\
09090-400 Santo Andr\'e - SP, Brazil\\ 
marcio.silva@ufabc.edu.br


\eject

\begin{figure}[ht]
\begin{center}
\includegraphics[width=12cm]{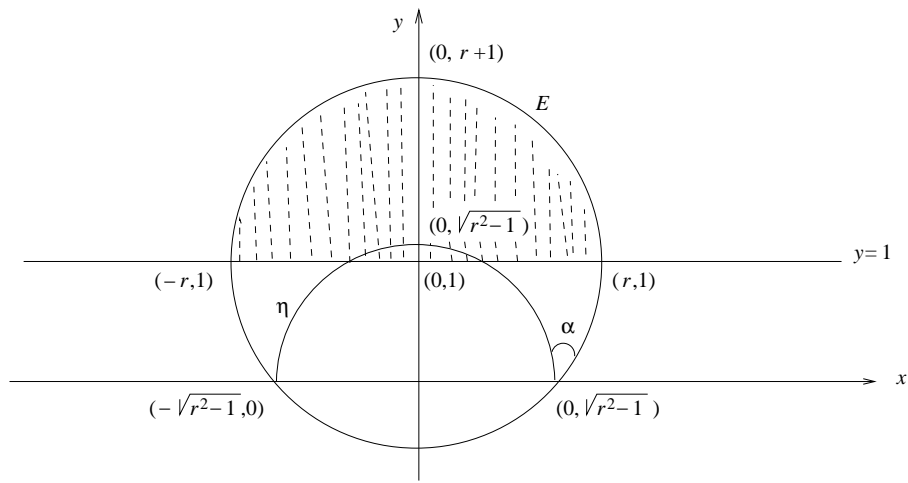}
\caption{Perimeter and area for an equidistant disk.}
\label{Fig:discoequidis}
\end{center}
\end{figure}

\end{document}